\documentclass[10pt]{article}
\usepackage{amssymb}
\usepackage{amsmath}
\usepackage{tikz}
\usepackage{longtable}
\usepackage{latexsym}
\usepackage{longtable}

\newcommand{\bea}{\begin{eqnarray}}
\newcommand{\eea}{\end{eqnarray}}
\newcommand{\be}{\begin{equation}}
\newcommand{\ee}{\end{equation}}
\newcommand{\bline}{\raisebox{2pt}{\tikz{\draw[-,black!40!black,solid,line width = 0.9pt](0,0) -- (7mm,0);}}}

\makeatletter\def\theequation{\arabic{section}.\arabic{equation}}\makeatother
\usepackage{pst-node}
\usepackage{pstricks}
\usepackage{graphicx}
\usepackage{latexsym}
\begin{document}

\title{\bf Gromov Product Decomposition of $7$-point Metric Spaces}
\author {{\bf  Ayse Humeyra Bilge\footnote{Faculty of Engineering and Natural Sciences, Department of Industrial Engineering, Kadir Has University, Istanbul, Turkey(ayse.bilge@khas.edu.tr)},
        \ \bf  Metehan Incegul\footnote{Istanbul Bilgi University, Faculty of Engineering and Natural Sciences, Kazim Karabekir Cad. No: 2/13  Eyup 34060 Istanbul / Turkey}
}}
          \maketitle

 \begin{abstract}
 \baselineskip=12pt \noindent
  Let $X$ be a finite metric space with elements $P_i$, $i=1,\dots ,n$ and with  distance functions $d_{ij}$. The Gromov product of the  triangle with vertices $P_i$, $P_j$ and $P_k$ at the vertex $P_i$ is defined by  $\Delta_{ijk}=\frac{1}{2}(d_{ij}+d_{ik}-d_{jk})$. A metric space is called $\Delta$-generic, if the set of   Gromov products at each $P_i$ has a unique smallest element $\Delta_{ijk}$.
  For a $\Delta$-generic metric space, the map $P_i\to (P_jP_k)$, where $(P_jP_k)$ is the edge joining $Pj$ to $P_k$ is a well defined map called the ``Gromov product structure" [Bilge, Celik and Kocak, ``An equivalence class decomposition of finite metric spaces", {\it Discrete Mathmetics}, Vol 340, (2017) 1928-1932].  For $n=5$, the   $3$   $\Delta$-equivalence classes coincide with the classification of $5$-point metrics. For $n=6$, there are $26$  $\Delta$-equivalence classes refined by $339$ metric classes. In this paper, we present the first systematic treatment of $7$-point spaces and we obtain the   $\Delta$-equivalence  decomposition of $7$-point metric spaces that consist of  $431$ equivalence classes. 
  \end{abstract}


\noindent\textit{Keywords\/}: {\ Finite metric spaces, Gromov product structure, Weighted graphs}

\baselineskip=16pt

\section{Introduction}

An $n$-point metric space  is described by the set of distances $d_{ij}$ subject to triangle
inequalities.  The structure of a finite metric is described by the so-called ``metric fan" and the metrics
belonging to the interior of a metric fan are called   ``generic" \cite{Sturmfels}."
In a series of papers,    we studied metric spaces in terms of the ``Gromov products" as defined by Eqn.(2.1) \cite{BCK2015},\cite{BCK2017},\cite{BCK2018},\cite{BI_Malta},\cite{IB_algorithm}. 
In this approach, distance functions are expressed in terms of Gromov products   and triangle inequalities are replaced by  the compatibility conditions of the defining equations for the $d_{ij}$'s as given by (2.2).  

As it is well known  \cite{Koolen}, \cite{Celik2018}, there are  $1$ and $3$ inequivalent metrics on  $4$-point and $5$-point spaces, respectively and the
$\Delta$-equivalence decomposition coincides with this classification.  For $n=6$,  there are $339$ inequivalent metrics
as classified in \cite{Sturmfels}. 
 We presented the $\Delta$-equivalence classes of $6$-point metrics in  \cite{BCK2017},\cite{IB_algorithm} where we have shown that
there are $26$ equivalence classes. 

 In the present work, we study the $\Delta$-equivalence decomposition of $7$-point spaces by using the algorithm of \cite{IB_algorithm} that   
consists of $4$ steps. At first, the algorithm enumerates all
Gromov product structures, then at step two, eliminates the ones that are not allowable using Proposition 1 and Corollary 1. The third step that consists of the application of the full permutation  group to allowable Gromov products is the most time consuming part of the algorithm.   
Finally at the fourth step, non-generic  Gromov product structures are eliminated \cite{IB_algorithm}.

 For $6$-point spaces
this algorithm runs on a standard computer in a reasonable time without the need for more sophisticated methods.
For $n=7$,  the third step of the algorithm, consisting of applying the full permutation group on $7$ points to all allowable
Gromov product structures turned out to be practically impossible to run on a standard computer.  In the search of permutation invariants of Gromov product
 structures, we ended up with a matrix representation, as presented in \cite{BI_Malta}.

In Section 2, we briefly review the basics of $\Delta$-equivalence decomposition.  
  The matrix representation of a Gromov product structure is defined in Section 3  and the invariants that are  used 
for the classification are described in Proposition 2. In Section 4, we illustrate the method for $6$-point spaces and in Section 5, we present  the 
classification of $7$-point spaces that consist of 
 $431$ Gromov product equivalence classes.  The complete list of these classes are presented in the Appendix.

\section{Preliminaries }
\def\theequation{2.\arabic{equation}}\makeatother
\setcounter{equation}{0}

Let $X$ be a finite metric space with $n$ elements $P_i$, $i=1,\dots,n$ and let $d_{ij}$ be the distance
between $P_i$ and $P_j$. We denote the triangle with vertices $P_i$, $P_j$ and $P_k$ by $(P_iP_jP_k)$ .
The quantity $\Delta_{ijk}$ defined by 
$$ \Delta_{ijk}=\Delta_{ikj}={\textstyle\frac{1}{2}}(d_{ij}+d_{ik}-d_{jk})\eqno(2.1)$$
is called the ``Gromov product" of the triangle $(P_iP_jP_k)$ at the
vertex $P_i.$

The triangle inequalities are equivalent to the
non-negativity of the $\Delta_{ijk}$'s. The solution set of the triangle inequalities is a convex polyhedral cone
in $R^{n(n-1)/2}$ and the classification of metric spaces is given in terms of the combinatorial
properties of this cone \cite{Sturmfels}.

For  any triangle $(P_iP_jP_k)$, the distances
can be expressed in terms of the Gromov products as
$$d_{ij}=\Delta_{ijk}+\Delta_{jik},\quad k=1,\dots n, \quad k\ne i,j. \eqno(2.2)$$
For an $n$-point  metric space, there are $n(n-1)(n-2)/2$ Gromov products that has to satisfy $n-3$ equalities for each
unordered pair $(i,j)$.  It follows that  an $n$-point
metric space is characterized by positive solutions of $n(n-1)(n-3)/2 $ linear equations in $n(n-1)(n-2)/2$ unknowns.
A finite metric space is called ``$\Delta$-generic", if the set of
Gromov products at each $P_i$ has a unique smallest element.
We define Gromov product structures on $\Delta$-generic spaces as follows.
\vskip 0.2cm\noindent
{\bf Definition 1.}{\it
Let $X$ be a $\Delta$-generic $n$-point  metric space with elements $P_a$, let $E(X)$ be the set of edges of $X$ and let
  $\Delta_{a,b_a,c_a}$ be the minimal Gromov Product at $P_a$. The Gromov product structure on $X$ is the map ${\cal E}$ from $X$ to the set $E(X)$
  defined by  ${\cal E}(P_a)=(P_{b_a}P_{c_a})$, where $(P_{a}P_{b})$ is the edge joining $P_a$ to $P_b$. 
 Two $\Delta$-generic  finite metric spaces are ``$\Delta$-equivalent" if they have the same Gromov product structures, up to permutations of indices.}
\vskip 0.2cm

The Gromov products at different $P_i$'s are not independent.  Their relations are given by the Proposition 1
below, whose proof is straightforward.
\vskip 0.2cm\noindent
{\bf Proposition 1.}{\it
Let  $\Delta_{abc}$ be the  Gromov product at node $a$ and $i$ be
a node different from $a$, $b$ and $c$.  Then, }
$$
\Delta_{abi}-\Delta_{abc}=\Delta_{cbi}-\Delta_{cai}
=\Delta_{iac}-\Delta_{ibc}=\Delta_{bac}-\Delta_{bai},\eqno(2.3a)$$
$$\Delta_{aci}-\Delta_{abc}=\Delta_{bci}-\Delta_{bai}
=\Delta_{iab}-\Delta_{ibc}=\Delta_{cab}-\Delta_{cai}.\eqno(2.3b)$$
\vskip 0.2cm
\noindent
These equalities lead to a convenient algorithmic procedure for determining allowable Gromov product structures.

\vskip 0.2cm\noindent
{\bf Corollary 1.}{\it
If $\Delta_{abc}$ is the minimal Gromov Product at
$P_a$ then,  $\Delta_{cbi}$ and $\Delta_{cab}$ cannot be
minimal at  $P_c$, $\Delta_{bci}$ and $\Delta_{bac}$ cannot be
minimal at  $P_b$ and $\Delta_{iac}$ and $\Delta_{iab}$ cannot
be minimal at  $P_i$.}
\vskip 0.2 cm
\noindent
  Thus, the minimality of $\Delta_{abc}$ leads to the exclusion of Gromov product structures containing
  $\Delta_{bci}$, $\Delta_{cbi}$, $\Delta_{iab}$, $\Delta_{iac}$ from the set of allowable Gromov product structures.
These elimination rules  lead to an algorithmic method for obtaining Gromov Product equivalence classes.
 We start by all possible Gromov Product index sets, apply the elimination rules and collect together the ones
 that are mapped to each other under permutations. This process works well for $n=5$ and $n=6$ but it becomes practically
 impossible for $n=7$.  For $n\ge 7$ we will eliminate beforehand a large number of equivalent cases.
 This will be done by using ``excess cycle lengths", to be defined below.

\vskip0.2cm\noindent
{\bf Definition 2.} {\it
If the Gromov products at nodes $\{a_1,a_2,\dots,a_k\}$ are
$$\{\Delta_{{a_1}{a_i}{a_2}},\Delta_{{a_2}{a_1}{a_3}}, \dots, \Delta_{{a_k},{a_{k-1}},{a_j}}\},$$
then we say that the Gromov Product structure
$$\{(a_1a_ia_2),(a_2a_1a_3), \dots, (a_k,a_{k-1},a_j)\}$$
contains a $\Delta$-chain of length $k$. If $k=1$, then the chain has length $1$ and $P_{a_1}$ is called an `isolated" point. For $k>1$ The points $P_{a_1}$ and $P_{a_k}$ are ``end points" of the chain while  the points $P_{a_i}$, $i\ne 1,k$ are ``interior points" of the chain.  For $k>1$, if $i=k$ and $j=1$, then we say that there is a $\Delta$-cycle of length $k$. A chain of length $2$ has no interior points.  A cycle has no end points. 
}
\vskip 0.2 cm
\noindent
The cycle and chain structure allows to choose Gromov product structures in certain canonical forms and reduces considerably the burden of finding equivalence classes.

\section {Matrix representation of Gromov product structures}

An $n$-point  metric space is represented by a weighted complete graph.  If $\Delta_{aij}$ is the minimal Gromov product as $P_a$,
then the edges joining $P_a$ to the other $P_b$'s can be sticked together by an amount $\Delta_{aij}$ and the edge joining $P_i$ to $P_j$ can be removed \cite{BCK2015}.
Repeating this procedure at each point, we obtain the so-called ``pendant-free" reduction of the complete graph.

In this section, we define the matrix representation of a Gromov product structure as defined in \cite{BI_Malta},
discuss its relation with the graph of pendant free reductions and we give the matrix representations for $n=4$ and $n=5$.

\vskip 0.2cm
\noindent{\bf Definition 3.}
{\it Let
$S=\{\Delta_{1a_1b_1},\Delta_{2a_2b_2},\dots,\Delta_{na_nb_n}\}$
be a Gromov product structure for a finite metric space with $n$ elements.
The  ``matrix representation for $S$" is the matrix $G_S$, defined by
$G_S(i,j)=1$, if the $j=a_i$ or $j=b_i$, in $\Delta_{ia_ib_i}$, and
$G_S(i,j)=0$ otherwise.}
\vskip 0.2cm
\noindent
Note that $G_S$ consists of zeros and ones only and it
has exactly two $1$'s in each row. It follows that $\lambda=2$ is always an eigenvalue with corresponding eigenvectors
multiples of $X=(1,1,\dots,1)^t$.
We note that similarity of the matrix representations don't imply equivalence of Gromov product structures, because the change of basis matrix need not be a permutation matrix. Nevertheless the (integer valued) quantities ${\rm Trace}(G_S^k)$ will be useful to discriminate between inequivalent Gromov product structures. 
\vskip 0.2cm
\noindent{\bf $4$-point spaces.}
For  $n=4$, there is unique Gromov product structure that can be given as
$$X^{(4)}: \quad \{\Delta_{124},\Delta_{213},\Delta_{324},\Delta_{413}\}$$
and its matrix is
$$G=\left(
\begin{array}{cccc}
0 & 1 & 0 & 1 \\
1 & 0 & 1 & 0 \\
0 & 1 & 0 & 1 \\
1 & 0 & 1 & 0 \\
\end{array}
\right)
$$
\vskip 0.2cm
\noindent{\bf $5$-point spaces.}
For $n=5$, there are $3$ Gromov product structures  corresponding  to the
metric classes
\begin{align*}
X^{(5)}_A&\quad\{\Delta_{125},\Delta_{213},\Delta_{324},\Delta_{435},\Delta_{514}\},\\
X^{(5)}_B&\quad \{\Delta_{124},\Delta_{213},\Delta_{324},\Delta_{435},\Delta_{524}\},\\
X^{(5)}_C&\quad \{\Delta_{124},\Delta_{213},\Delta_{324},\Delta_{413},\Delta_{513}\}.
\end{align*}
and their matrix representations are 
\begin{align*}
G_{A}=\left(
\begin{array}{ccccc}
0 & 1 & 0 & 0 & 1\\
1 & 0 & 1 & 0 & 0\\
0 & 1 & 0 & 1 & 0\\
0 & 0 & 1 & 0 & 1\\
1 & 0 & 0 & 1 & 0\\
\end{array}        \right),\quad 
&G_{B}=\left(
\begin{array}{ccccc}
0 & 1 & 0 & 1 & 0\\
1 & 0 & 1 & 0 & 0\\
0 & 1 & 0 & 1 & 0\\
0 & 0 & 1 & 0 & 1\\
0 & 1 & 0 & 1 & 0\\
\end{array}         \right),\quad
&G_{C}=\left(
\begin{array}{ccccc}
0 & 1 & 0 & 1 & 0\\
1 & 0 & 1 & 0 & 0\\
0 & 1 & 0 & 1 & 0\\
1 & 0 & 1 & 0 & 0\\
1 & 0 & 1 & 0 & 0\\
\end{array}          \right).
\end{align*}
Their  $\Delta$-chain and $\Delta$-cycle structures are given by
\begin{align*}
X^{(5)}_A&\quad
\stackrel{5}{\circ}\bline\stackrel{1}{\bullet}\bline\stackrel{2}{\bullet}\bline\stackrel{3}{\bullet}\bline\stackrel{4}{\bullet}\bline\stackrel{5}{\bullet}\bline\stackrel{1}{\circ}\\
X^{(5)}_B&\quad
\stackrel{4}{\circ}\bline\stackrel{1}{\bullet}\bline\stackrel{2}{\bullet}\bline\stackrel{3}{\bullet}\bline\stackrel{4}{\bullet}\bline\stackrel{5}{\bullet}\bline\stackrel{2}{\circ}\\
X^{(5)}_C&\quad
\stackrel{4}{\circ}\bline\stackrel{1}{\bullet}\bline\stackrel{2}{\bullet}\bline\stackrel{3}{\bullet}\bline\stackrel{4}{\bullet}\bline\stackrel{1}{\circ}  \quad \quad
\stackrel{1}{\circ}\bline\stackrel{5}{\bullet}\bline\stackrel{3}{\circ}.
\end{align*}

We will now present certain results on the relation of matrix invariants to the properties of pendant free reductions.

An $n\times n$  matrix $A$ is irreducible, if $\sum_{i=0}^n A^i$ has no zero element.   

\vskip 0.2cm\noindent
{\bf Definition 4.}{\it
A Gromov product structure is called ``irreducible" if the corresponding matrix is irreducible. 
Otherwise it is called ``reducible".}
\vskip 0.2cm\noindent
If a metric space is reducible in the sense above, then,  a proper subset is itself a metric space.  Hence, the parametrization of this subspace can be enlarged to the parametrization of the larger space.

\vskip 0.2cm\noindent
{\bf Proposition 2.}{\it
Let $G_S$ be the matrix representation of a Gromov product structure $S$, as defined in Definition 3. Then,
\begin{description}
\item{(i)} The rank of $G_S$ is the number of edges removed in the pendant free reduction.
\item{(ii)} If $G_S$ is reducible, then the number of columns with no zero element  in $\sum_{i=0}^n G_S^i$ 
gives the dimension of the largest irreducible subspace.    
\item{(iii)} The symmetric part of the matrix $G_S$, expressed as $H_S={\rm Floor}(\frac{1}{2}(G_S+G_S^t))$ gives information of the chain structure. The sum of the entries of the columns are $0$, $1$ or $2$.  $P_i$ is respectively an ``isolated point", an ``end point" or an ``interior point" according as the sum of the entries of the $i$'th column of $H_S$ is $0$, $1$ or $2$. 
\end{description}
}
\vskip 0.2cm
\noindent

\section{Gromov Product Classification for $n=6$}
In previous work, we obtained the $\Delta$-equivalence classification of $6$-point metrics by directly applying the algorithm described in Section 1 \cite{BI_Malta}. 
In this section, we reiterate this procedure by using the notion of  ``cycle" and ``chain" lengths,  in order to display their  advantages.
The structure of $\Delta$-chains and  $\Delta$-cycles  (Definition 2) on an $n$-point space can be conveniently analysed by the ``partition of integers".  
The partition of integers for $n=6$ is
\begin{align*}
6+0&& 5+1 && 4+1+1&&  3+1+1+1 && 2+1+1+1+1&& 1\times 6  \\
 && 4+2 && 3+2+1&&  2+2+1+1 &&          &&               \\
 && 3+3 && 2+2+2&&          &&          &&
 \end{align*}
If there is a $6$-cycle of $6$-chain, without loss of generality, the chain structure can be chosen as
$$
	 \stackrel{i}{\circ}\bline\stackrel{1}{\bullet}\bline\stackrel{2}{\bullet}\bline\stackrel{3}{\bullet}\bline\stackrel{4}{\bullet}
\bline\stackrel{5}{\bullet}\bline\stackrel{6}{\bullet}\bline\stackrel{j}{\circ}$$
Here, if $i=6$ and $j=1$ there is a $6$-cycle; otherwise there is a $6$-chain. If there is a cycle, all points are interior points. In the case of a chain, there are $2$ end points and $4$ interior points.

If there is a $5$-cycle or a  $5$-chain, the chain structure is
$$
 \stackrel{i}{\circ}\bline\stackrel{1}{\bullet}\bline\stackrel{2}{\bullet}\bline\stackrel{3}{\bullet}\bline\stackrel{4}{\bullet}
\bline\stackrel{5}{\bullet}\bline\stackrel{j}{\circ}\quad \quad   \stackrel{k}{\circ}\bline\stackrel{6}{\bullet}\bline\stackrel{l}{\circ} $$
Here, if there is a $5$-cycle, we have $5$ interior points and $1$ isolated point. If there is a $5$-chain, there are $2$ interior points, $2$ end points and $1$ isolated point.  	

For $n=6$, not all partitions of the integer $6$ leads to allowable types. 
The list of allowable types for $n=6$ are below. The numbers in the last column gives the number of isolated points, end points and  interior points respectively. Note that $4+2$ (Chain) and $3+3$ types are not distinguished by these invariants. 
\newpage
\begin{table}[]
	\centering
	\caption{Allowable cycle and chain configurations for $6$-point spaces}
	\label{my-label}
	\begin{tabular}{|l|l|l|}
		\hline
$6+0$(cycle) &$\quad\circ\bline\bullet\bline\bullet\bline\bullet\bline\bullet\bline\bullet\bline\bullet\bline\circ$
&(0,0,6)\\
\hline
$6+0$(chain) &$\quad\circ\bline\bullet\bline\bullet\bline\bullet\bline\bullet\bline\bullet\bline\bullet\bline\circ$
&(0,2,4)\\
\hline
$5+1$(cycle)	&$\quad\circ\bline\bullet\bline\bullet\bline\bullet\bline\bullet\bline\bullet\bline\circ$&(1,0,5)\\
&$\quad\circ\bline\bullet\bline\circ$&\\
\hline
$5+1$(chain)	&$\quad\circ\bline\bullet\bline\bullet\bline\bullet\bline\bullet\bline\bullet\bline\circ$&(1,2,3)\\
&$\quad\circ\bline\bullet\bline\circ$&\\
\hline
$4+2$(cycle)	&$\quad\circ\bline\bullet\bline\bullet\bline\bullet\bline\bullet\bline\circ$&(0,2,4)\\
&$\quad\circ\bline\bullet\bline\bullet\bline\circ$&\\
\hline
$4+2$(chain)	&$\quad\circ\bline\bullet\bline\bullet\bline\bullet\bline\bullet\bline\circ$&(0,4,2)\\
&$\quad\circ\bline\bullet\bline\bullet\bline\circ$&\\
\hline
$3+3$&$\quad\circ\bline\bullet\bline\bullet\bline\bullet\bline\circ$&(0,4,2)\\
&$\quad\circ\bline\bullet\bline\bullet\bline\bullet\bline\circ$&\\
\hline
$4+1+1$(cycle)&$\quad\circ\bline\bullet\bline\bullet\bline\bullet\bline\bullet\bline\circ$&(2,0,4)\\
&$\quad\circ\bline\bullet\bline\circ\quad\quad\circ\bline\bullet\bline\circ$&\\
\hline
$4+1+1$(chain)&$\quad\circ\bline\bullet\bline\bullet\bline\bullet\bline\bullet\bline\circ$&(2,2,2)\\
&$\quad\circ\bline\bullet\bline\circ\quad\quad\circ\bline\bullet\bline\circ$&\\
\hline
$3+2+1$&$\quad\circ\bline\bullet\bline\bullet\bline\bullet\bline\bullet\bline\circ$&(1,4,1)\\
&$\quad\circ\bline\bullet\bline\bullet\bline\circ\quad\quad\circ\bline\bullet\bline\circ$&\\
\hline
$2+2+2$&$\quad\circ\bline\bullet\bline\bullet\bline\circ$&(0,6,0)\\
&$\quad\circ\bline\bullet\bline\bullet\bline\circ\quad\quad\circ\bline\bullet\bline\bullet\bline\circ$&\\
\hline
$3+1\times 3$&$\quad\circ\bline\bullet\bline\bullet\bline\bullet\bline\circ$&(3,1,2)\\
&$\quad\circ\bline\bullet\bline\circ\quad\quad\circ\bline\bullet\bline\circ\quad\quad\circ
\bline\bullet\bline\circ$&\\
\hline
$2+2+1+1$&$\quad\circ\bline\bullet\bline\bullet\bline\circ
\quad\quad\circ\bline\bullet\bline\bullet\bline\circ$&(2,4,0)\\
&$\quad\circ\bline\bullet\bline\circ\quad\quad\circ\bline\bullet\bline\circ$&\\
\hline
$2+1\times 4$&$\quad\circ\bline\bullet\bline\bullet\bline\circ$&(4,2,0)\\
&$\quad\circ\bline\bullet\bline\circ\quad\quad\circ\bline\bullet\bline\circ$&\\
&$\quad\circ\bline\bullet\bline\circ\quad\quad\circ\bline\bullet\bline\circ$&\\
\hline
\end{tabular}
\end{table}
\newpage
The complete list of $6$-point spaces grouped according to their cycle and chain structures is given below. In this table we use $abc$ to denote the minimal Gromov product $\Delta_{abc}$ at $P_a$. 

\newpage
\begin{table}[h]
	\centering
	\caption{Complete list of $6$-point spaces}
	\label{my-label}
	\begin{tabular}{|l|l|l|l|l|l|l|l|}
		\hline
		{\rotatebox{0}{Type}}&
		{\rotatebox{90}{R/I}}&
	    {\rotatebox{90}{Int.}}&
	    {\rotatebox{90}{End.}}&
	    {\rotatebox{90}{Disc.}}&
	    {\rotatebox{90}{Remv.}}&  
        {\rotatebox{0}{Gromov Product Structure}}&   Type                                                              \\
\hline
6+0 (Cycle)				& I & 6& 0& 0   & 6 &126,213,324,435,546,615&$(I_{17})$ \\
\hline
6+0 (Chain) 		    & I & 4& 2& 0   & 4 &124,213,324,435,546,635&$I_{7}$ \\
  	                    & I &  &  &     & 5 &124,213,324,435,546,615&$I_{8}$ \\
						& I &  &  &     & 5 &124,213,324,435,546,625&$I_{14}$ \\
						& I &  &  &     & 5 &125,213,324,435,546,625&$I_{15}$\\
\hline
5+1 (Cycle)  			& R & 5& 0& 1   & 5 &125,213,324,435,514,613&$R_{8}$ \\ 
\hline
5+1 (Chain)   			& I & 3& 2& 1  & 5 &125,213,324,435,546,613 &$I_{11}$ \\
                        & I &  &  &     & 4 &124,213,324,435,546,613&$I_{13}$ \\
						& R &  &  &     &3 &124,213,324,435,524,624 &$R_{2}$\\
						& R &  &  &     &3 &124,213,324,435,524,613 &$R_{4}$ \\
						& R &  &  &     &4 &124,213,324,435,524,615 &$R_{6}$ \\ 
\hline
4+2 (Cycle) 			& R & 4& 2& 0  & 4 &124,213,324,413,516,625 &$R_7$ \\
     					& R &  &  &    & 4 &124,213,324,413,516,635 &$R_9$ \\ 
\hline
4+2 (Chain) 			& I & 2& 4& 0   &4 &124,213,324,435,526,635 &$I_5$ \\
                        & I &  &  &     &4 &125,213,324,436,536,625 &$I_6$ \\
						& I &  &  &    & 4 &124,213,324,435,516,635 &$I_9$ \\
						& I &  &  &     &5 &124,213,324,435,516,625 &$I_{10}$ \\
						& I &  &  &    & 5 &124,213,324,435,526,615 &$I_{12}$ \\ 
\hline
4+1+1 (Cycle)			& R & 4& 0& 2  & 2 &124,213,324,413,513,613 &$R_1$ \\
    					& R &  &  &    & 2 &124,213,324,413,513,624 &$R_3$ \\
						& R &  &  &    & 3 &124,213,324,413,513,625 &$R_5$ \\ 
\hline
4+1+1 (Chain) 			& I & 2& 2& 2  & 4 &124,213,324,435,516,624 &$I_3$ \\ 
\hline
3+3        				& I & 2& 4& 0  & 3 &124,213,324,456,524,624 &$I_1$ \\ 
\hline
3+2+1       			& I & 1& 4& 1  & 3&124,213,324,456,513,624 &$I_2$  \\ 
\hline
3+1+1+1    				& I & 1& 2& 3  & 3&124,213,324,456,513,613 &$I_4$ \\ 
\hline
2+2+2      		    	& I & 0& 6& 0  & 3 &156,213,324,456,513,624 &$I_{16}$ \\ 
\hline
	\end{tabular}
\end{table}
\newpage
\section{Gromov Product Chain Classification for $n=7$}

For $n=7$, the partition of integers lead to the following chain configurations.

The partition of integers for $n=7$ is
\begin{align*}
7+0&& 6+1 && 5+1+1&&  4+1+1+1 && 3+1+1+1+1&& 2+1\times 5 &&1\times 7  \\
 && 5+2 && 4+2+1&&  3+2+1+1 && 2+2+1+1+1&&                          \\
 && 4+3 && 3+3+1&&  2+2+2+1 &&          &&
\end{align*}

Without loss of generality, we choose the following canonical orderings.  Algorithmically it is possible to terminate the chains. But,  even if we start with canonical forms with shorter chains, it is difficult to avoid the occurrence of longer chains.  The grouping of the Gromov products by  their matrix invariants  help with this problem.

\noindent There is a $7$ cycle or a $7$ chain: 
$$
\stackrel{i}{\circ}\bline\stackrel{1}{\bullet}\bline\stackrel{2}{\bullet}\bline\stackrel{3}{\bullet}\bline\stackrel{4}{\bullet}
\bline\stackrel{5}{\bullet}\bline\stackrel{6}{\bullet}\bline\stackrel{7}{\bullet}\bline\stackrel{j}{\circ}$$
	There is a chain of length $6$ which is not extendible to a longer chain : 
$$
\stackrel{i}{\circ}\bline\stackrel{1}{\bullet}\bline\stackrel{2}{\bullet}\bline\stackrel{3}{\bullet}\bline\stackrel{4}{\bullet}
\bline\stackrel{5}{\bullet}\bline\stackrel{6}{\bullet}\bline\stackrel{j}{\circ}$$	
 There is a chain of length $5$ which is not extendible to a longer chain : 
$$
\stackrel{i}{\circ}\bline\stackrel{1}{\bullet}\bline\stackrel{2}{\bullet}\bline\stackrel{3}{\bullet}\bline\stackrel{4}{\bullet}
\bline\stackrel{5}{\bullet}\bline\stackrel{j}{\circ}$$	
There is a chain of length $4$ which is not extendible to a longer chain : 
	$$
	\stackrel{i}{\circ}\bline\stackrel{1}{\bullet}\bline\stackrel{2}{\bullet}\bline\stackrel{3}{\bullet}\bline\stackrel{4}{\bullet}\bline\stackrel{j}{\circ}$$
There is a chain of length $3$ which is not extendible to a longer chain : 
	$$
	\stackrel{i}{\circ}\bline\stackrel{1}{\bullet}\bline\stackrel{2}{\bullet}\bline\stackrel{3}{\bullet}\bline\stackrel{4}{\circ}$$
There is a chain of length $2$ which is not extendible to a longer chain : 
	$$
	\stackrel{1}{\circ}\bline\stackrel{2}{\bullet}\bline\stackrel{3}{\bullet}\bline\stackrel{4}{\circ}$$
The longest chain of length $2$: 
		$$
		\stackrel{1}{\circ}\bline\stackrel{2}{\bullet}\bline\stackrel{3}{\bullet}\quad
		\stackrel{4}{\circ}\bline\stackrel{3}{\bullet}\bline\stackrel{5}{\bullet}
	$$     	
	The $\Delta$-equivalence decomposition of $7$-point spaces is obtained by implementing the algorithms described below with MATLAB. Details of the algorithms  are given in \cite{IB_algorithm}. We adapted the codes to $n=7$. 
For each chain length,  we start with canonical types generate all possible Gromov product structures. Then use Proposition 1 to select allowable ones.  At this stage, the sets of allowable structures are too big for applying the full permutation group.  Therefore, we use matrix invariants to obtain smaller groups to which permutation group acts in a reasonable time.  
 This way the problem is reduced to a feasible size.  Finally, we eliminate the structures that are not generic.  As a result we obtain $431$ types listed in the Appendix. The number of types for each chain length is listed below.

\begin{table}[htp]
	\centering
	\caption{$7$ point spaces}
	\label{my-label}
	\begin{tabular}{|l|l|l|l|l|l|}
		\hline
		{\rotatebox{0}{Type}}&
		{\rotatebox{90}{R/I}}&
		{\rotatebox{90}{Isolated}}&
		{\rotatebox{90}{End}}&
		{\rotatebox{90}{Internal}}&
		{\rotatebox{90}{Number}}	
		\\
		\hline
		7 (Cycle)  & I  & 0& 0& 7 & 1 \\ 		\hline
		7 (Chain)   & I & 0& 2& 5    & 8 \\		
		\hline
		6+1 & I& 1& 2& 4    & 8 \\     
		\hline		      
		5+2       & I   & 0& 4& 3                               &43 \\ 
		4+3        & I  & 0& 4& 3                               &   \\        \hline
		5+1+1     & I   & 2& 2& 3                               &22 \\		\hline
		4+2+1       & I & 1& 4& 2                               &63 \\
		3+3+1       & I & 1& 4& 2                               &   \\        \hline
		3+2+2       & I & 0& 6& 1                               &27 \\        \hline
		4+1+1+1 & I & 3& 2& 2                               &27 \\        \hline
		3+2+1+1    & I  & 2& 4& 1                               &42 \\        \hline
		2+2+2+1     & I & 1& 6& 0                               &14 \\        \hline
		3+1+1+1+1 & I & 4& 2& 1                               &12 \\        \hline
		2+2+1+1+1 & I& 3& 4& 0                               &13 \\        \hline
		2+1+1+1+1+1 & I& 5& 2& 0                               &6  \\        \hline
		1$\times$     7& I& 7& 9& 0                               &1  \\        \hline
	
		Contains $X^{(4)}$ & R & & & &  24 \\ 		\hline
		Contains $X^{(5)}$ & R & & & &  27 \\ 		\hline
		Contains $X^{(6)}$ & R & & & &  93 \\ 		\hline
	\end{tabular}
\end{table}

\newpage
\appendix
\setcounter{table}{0}
\renewcommand{\thetable}{A\arabic{table}}
\section{List of 7-Spaces}

\begin{table}[htp]\centering\caption{7-Cycle and 7-Chains}\begin{tabular}{|l|l|l|l|l|l|l|}\hline 	 
127&   213&   324&   435&   546&   657&   716\\ 
\hline 
124 &  213 &  324  & 435 &  546 &  657 &  716\\		
125 &  213 &  324  & 435 &  546 &  657 &  716\\
124 &  213 &  324  & 435 &  546 &  657 &  736\\
124 &  213 &  324  & 435 &  546 &  657 &  726\\
126 &  213 &  324  & 435 &  546 &  657 &  726\\
125 &  213 &  324  & 435 &  546 &  657 &  726\\
125 &  213 &  324  & 435 &  546 &  657 &  736\\
124 &  213 &  324  & 435 &  546 &  657 &  746\\
\hline	\end{tabular}\end{table}
\begin{table}[htp]\centering\caption{6+1 Decomposition}\begin{tabular}{|l|l|l|l|l|l|l|}\hline 
124 &  213 &  324  & 435 &  546 &  657 &  713\\     
124 &  213 &  324  & 435 &  546 &  657 &  724\\
125 &  213 &  324  & 435 &  546 &  657 &  713\\
125 &  213 &  324  & 435 &  546 &  657 &  724\\
125 &  213 &  324  & 435 &  546 &  657 &  714\\
126 &  213 &  324  & 435 &  546 &  657 &  713\\
126 &  213 &  324  & 435 &  546 &  657 &  714\\
126 &  213 &  324  & 435 &  546 &  657 &  724\\		
\hline	\end{tabular}\end{table}
\begin{table}[htp]\centering\caption{5+2 Decomposition}\begin{tabular}{|l|l|l|l|l|l|l|}\hline 		
	124&	213&	324&	435&	546&	627&	716\\
	124&	213&	324&	435&	546&	637&	716\\
	124&	213&	324&	435&	546&	637&	726\\
	124&	213&	324&	435&	547&	627&	716\\
	124&	213&	324&	435&	547&	637&	716\\
	124&	213&	324&	435&	547&	637&	726\\
	124&	213&	324&	435&	547&	647&	716\\
	124&	213&	324&	435&	547&	647&	726\\
	124&	213&	324&	435&	547&	647&	736\\
	125&	213&	324&	435&	546&	627&	716\\
	125&	213&	324&	435&	546&	637&	716\\
	125&	213&	324&	435&	546&	637&	726\\
	125&	213&	324&	435&	547&	627&	716\\
	125&	213&	324&	435&	547&	637&	716\\
	125&	213&	324&	435&	547&	637&	726\\
	125&	213&	324&	435&	547&	647&	716\\
	125&	213&	324&	435&	547&	647&	726\\
	125&	213&	324&	435&	547&	647&	736\\
	126&	213&	324&	435&	546&	637&	726\\
	126&	213&	324&	435&	547&	637&	726\\
	126&	213&	324&	435&	547&	647&	726\\
\hline	\end{tabular}\end{table}
\begin{table}[htp]\centering\caption{4+3 Decomposition}\begin{tabular}{|l|l|l|l|l|l|l|}\hline 		
	124&	213&	324&	435&	567&	615&	715\\
	124&	213&	324&	435&	567&	625&	715\\
	124&	213&	324&	435&	567&	625&	725\\
	124&	213&	324&	435&	567&	635&	715\\
	124&	213&	324&	435&	567&	635&	725\\
	124&	213&	324&	435&	567&	635&	735\\
	124&	213&	324&	436&	567&	615&	715\\
	124&	213&	324&	436&	567&	625&	715\\
	124&	213&	324&	436&	567&	625&	725\\
	124&	213&	324&	437&	567&	625&	715\\
	124&	213&	324&	437&	567&	635&	715\\
	124&	213&	324&	437&	567&	635&	725\\
	124&	213&	324&	437&	567&	645&	715\\
	124&	213&	324&	437&	567&	645&	725\\
	125&	213&	324&	435&	567&	625&	725\\
	125&	213&	324&	435&	567&	635&	725\\
	125&	213&	324&	436&	536&	657&	716\\
	125&	213&	324&	436&	567&	625&	725\\
	125&	213&	324&	437&	536&	657&	716\\
	125&	213&	324&	437&	536&	657&	726\\
	125&	213&	324&	437&	546&	657&	716\\
	125&	213&	324&	437&	567&	635&	725\\
\hline	\end{tabular}\end{table}
\begin{table}[htp]\centering\caption{5+1+1 Decomposition}\begin{tabular}{|l|l|l|l|l|l|l|}\hline 		
124&	213&	324&	435&	546&	627&	713\\
124&	213&	324&	435&	546&	627&	715\\
124&	213&	324&	435&	546&	637&	715\\
124&	213&	324&	435&	546&	637&	724\\
124&	213&	324&	435&	546&	637&	725\\
124&	213&	324&	435&	547&	624&	716\\
124&	213&	324&	435&	547&	625&	716\\
124&	213&	324&	435&	547&	635&	716\\
124&	213&	324&	435&	547&	635&	726\\
125&	213&	324&	435&	546&	627&	713\\
125&	213&	324&	435&	546&	637&	714\\
125&	213&	324&	435&	546&	637&	724\\
125&	213&	324&	435&	546&	637&	725\\
125&	213&	324&	435&	547&	624&	716\\
125&	213&	324&	435&	547&	625&	716\\
125&	213&	324&	435&	547&	635&	716\\
125&	213&	324&	435&	547&	635&	726\\
126&	213&	324&	435&	546&	637&	714\\
126&	213&	324&	435&	546&	637&	715\\
126&	213&	324&	435&	546&	637&	724\\
126&	213&	324&	435&	547&	635&	713\\
126&	213&	324&	435&	547&	647&	713\\
\hline	\end{tabular}\end{table}

\begin{table}[htp]\centering\caption{4+2+1 Decomposition}\begin{tabular}{|l|l|l|l|l|l|l|}\hline 		
124	&213	&324	&435	&567	&615	&713\\
124	&213	&324	&435	&567	&624	&715\\
124	&213	&324	&435	&567	&625	&713\\
124	&213	&324	&435	&567	&625	&724\\
124	&213	&324	&435	&567	&635	&713\\
124	&213	&324	&435	&567	&635	&724\\
124	&213	&324	&436	&527	&627	&715\\
124	&213	&324	&436	&537	&625	&715\\
124	&213	&324	&436	&537	&627	&715\\
124	&213	&324	&436	&547	&625	&715\\
124	&213	&324	&436	&547	&627	&715\\
124	&213	&324	&436	&567	&615	&713\\
124	&213	&324	&436	&567	&625	&713\\
124	&213	&324	&436	&567	&625	&724\\
124	&213	&324	&437	&526	&615	&715\\
124	&213	&324	&437	&526	&635	&715\\
124	&213	&324	&437	&526	&645	&715\\
124	&213	&324	&437	&536	&615	&715\\
124	&213	&324	&437	&536	&625	&715\\
124	&213	&324	&437	&536	&625	&725\\
124	&213	&324	&437	&536	&645	&715\\
124	&213	&324	&437	&536	&645	&725\\
124	&213	&324	&437	&546	&615	&715\\
124	&213	&324	&437	&546	&625	&715\\
124	&213	&324	&437	&546	&625	&725\\
124	&213	&324	&437	&546	&635	&715\\
124	&213	&324	&437	&546	&635	&725\\
124	&213	&324	&437	&567	&624	&715\\
125	&213	&324	&435	&567	&625	&713\\
125	&213	&324	&435	&567	&625	&714\\
125	&213	&324	&435	&567	&635	&713\\
125	&213	&324	&436	&536	&627	&716\\
125	&213	&324	&436	&536	&657	&713\\
125	&213	&324	&436	&537	&627	&716\\
125	&213	&324	&436	&547	&627	&716\\
125	&213	&324	&436	&567	&625	&713\\
125	&213	&324	&436	&567	&625	&714\\
125	&213	&324	&437	&536	&625	&716\\
125	&213	&324	&437	&536	&625	&725\\
125	&213	&324	&437	&536	&627	&716\\
125	&213	&324	&437	&537	&627	&716\\
125	&213	&324	&437	&537	&637	&716\\
125	&213	&324	&437	&546	&625	&716\\
125	&213	&324	&437	&546	&627	&716\\
\hline	\end{tabular}\end{table}
\newpage

\begin{table}[htp]\centering\caption{3+3+1 Decomposition}\begin{tabular}{|l|l|l|l|l|l|l|}\hline 		
124	&213	&324	&456	&524	&647	&715\\
124	&213	&324	&456	&547	&624	&713\\
124	&213	&324	&456	&547	&637	&715\\
124	&213	&324	&456	&547	&647	&713\\
124	&213	&324	&457	&546	&615	&713\\
124	&213	&324	&457	&546	&625	&713\\
124	&213	&324	&467	&526	&645	&715\\
124	&213	&324	&467	&567	&625	&715\\
124	&213	&324	&467	&567	&635	&715\\
124	&213	&324	&467	&567	&645	&713\\
125	&213	&324	&456	&537	&647	&716\\
125	&213	&324	&456	&547	&614	&713\\
125	&213	&324	&456	&547	&624	&713\\
125	&213	&324	&457	&546	&637	&714\\
125	&213	&324	&467	&536	&624	&714\\
125	&213	&324	&467	&536	&645	&713\\
125	&213	&324	&467	&537	&624	&714\\
125	&213	&324	&467	&537	&645	&714\\
125	&213	&324	&467	&537	&645	&724\\
\hline	\end{tabular}\end{table}

\begin{table}[htp]\centering\caption{3+2+2 Decomposition}\begin{tabular}{|l|l|l|l|l|l|l|}\hline 		
   124&   213&   324&   456&   527&   624&   715\\
   125&   213&   324&   417&   536&   625&   724\\
   125&   213&   324&   467&   536&   625&   724\\
   124&   213&   324&   467&   527&   624&   715\\
   124&   213&   324&   456&   537&   624&   715\\
   125&   213&   324&   467&   567&   625&   714\\
   125&   213&   324&   457&   536&   625&   724\\
   125&   213&   324&   457&   536&   657&   724\\
   124&   213&   324&   467&   567&   624&   715\\
   124&   213&   324&   467&   527&   645&   715\\
   124&   213&   324&   467&   537&   645&   715\\
   124&   213&   324&   456&   527&   647&   715\\
   124&   213&   324&   467&   567&   625&   724\\
   125&   213&   324&   467&   567&   625&   724\\
   125&   213&   324&   457&   536&   625&   714\\
   125&   213&   324&   417&   546&   657&   724\\
   125&   213&   324&   467&   536&   625&   714\\
   125&   213&   324&   417&   536&   657&   724\\
   125&   213&   324&   467&   567&   635&   714\\
   125&   213&   324&   457&   536&   657&   714\\
   125&   213&   324&   457&   546&   627&   716\\
   125&   213&   324&   456&   547&   627&   716\\
   125&   213&   324&   457&   546&   637&   716\\
   125&   213&   324&   456&   547&   637&   716\\
   125&   213&   324&   417&   536&   645&   724\\
   125&   213&   324&   417&   567&   645&   724\\
   125&   213&   324&   417&   536&   645&   745\\
\hline	\end{tabular}\end{table}
\newpage\newpage
\begin{table}[htp]\centering\caption{4+1+1+1 Decomposition}\begin{tabular}{|l|l|l|l|l|l|l|}\hline 		
   124&   213&   324&   435&   567&   613&   713\\
   124&   213&   324&   436&   527&   615&   713\\
   124&   213&   324&   436&   547&   615&   713\\
   124&   213&   324&   437&   526&   637&   715\\
   125&   213&   324&   436&   547&   625&   713\\
   125&   213&   324&   437&   537&   625&   716\\
   124&   213&   324&   435&   567&   624&   713\\
   124&   213&   324&   436&   547&   625&   713\\
   124&   213&   324&   437&   536&   624&   715\\
   124&   213&   324&   436&   524&   627&   715\\
   124&   213&   324&   436&   536&   627&   715\\
   124&   213&   324&   437&   536&   627&   715\\
   124&   213&   324&   437&   546&   627&   715\\
   125&   213&   324&   435&   567&   614&   713\\
   124&   213&   324&   435&   526&   637&   715\\
   125&   213&   324&   436&   537&   625&   714\\
   124&   213&   324&   437&   546&   637&   715\\
   124&   213&   324&   435&   526&   647&   715\\
   125&   213&   324&   436&   537&   625&   716\\
   124&   213&   324&   435&   567&   624&   724\\
   124&   213&   324&   437&   546&   637&   725\\
   125&   213&   324&   435&   567&   613&   713\\
   125&   213&   324&   435&   567&   624&   713\\
   125&   213&   324&   436&   547&   627&   713\\
   125&   213&   324&   435&   567&   614&   714\\
   125&   213&   324&   436&   537&   627&   714\\
   125&   213&   324&   436&   547&   625&   716\\
\hline	\end{tabular}\end{table}

\begin{table}[htp]\centering\caption{3+2+1+1 Decomposition}\begin{tabular}{|l|l|l|l|l|l|l|}\hline 		
   124&   213&   324&   457&   546&   613&   713\\
   124&   213&   324&   467&   567&   625&   713\\
   124&   213&   324&   467&   524&   624&   715\\
   125&   213&   324&   467&   567&   624&   713\\
   124&   213&   324&   457&   526&   615&   713\\
   125&   213&   324&   416&   537&   637&   716\\
   125&   213&   324&   457&   536&   624&   724\\
   124&   213&   324&   467&   526&   615&   713\\
   125&   213&   324&   467&   567&   614&   713\\
   124&   213&   324&   457&   536&   615&   713\\
   124&   213&   324&   467&   567&   615&   713\\
   124&   213&   324&   456&   527&   624&   713\\
   124&   213&   324&   457&   546&   627&   713\\
   124&   213&   324&   467&   526&   615&   715\\
   124&   213&   324&   456&   547&   627&   713\\
   125&   213&   324&   417&   536&   624&   724\\
   124&   213&   324&   467&   536&   615&   715\\
   125&   213&   324&   457&   546&   613&   713\\
   125&   213&   324&   416&   547&   624&   713\\
   125&   213&   324&   416&   547&   647&   713\\
   125&   213&   324&   416&   547&   624&   716\\
   124&   213&   324&   467&   536&   624&   715\\
   125&   213&   324&   416&   567&   624&   713\\
   125&   213&   324&   416&   537&   624&   716\\
   124&   213&   324&   467&   536&   625&   715\\
   124&   213&   324&   467&   526&   635&   715\\
   125&   213&   324&   467&   536&   614&   713\\
   125&   213&   324&   467&   536&   624&   713\\
   125&   213&   324&   457&   546&   627&   713\\
   125&   213&   324&   457&   536&   614&   714\\
   125&   213&   324&   457&   536&   624&   714\\
   125&   213&   324&   456&   537&   624&   714\\
   125&   213&   324&   457&   536&   627&   714\\
   125&   213&   324&   456&   537&   624&   716\\
   125&   213&   324&   457&   536&   625&   716\\
   125&   213&   324&   457&   536&   627&   716\\
   125&   213&   324&   416&   537&   627&   716\\
   125&   213&   324&   416&   547&   627&   716\\
   125&   213&   324&   456&   537&   637&   716\\
   125&   213&   324&   416&   537&   645&   716\\
   125&   213&   324&   416&   537&   645&   724\\
   125&   213&   324&   417&   536&   645&   726\\
\hline	\end{tabular}\end{table}

\begin{table}[htp]\centering\caption{2+2+2+1 Decomposition}\begin{tabular}{|l|l|l|l|l|l|l|}\hline 		
   145&   213&   324&   467&   513&   624&   713\\
   156&   213&   324&   457&   524&   613&   713\\
   145&   213&   324&   467&   516&   624&   713\\
   145&   213&   324&   416&   537&   637&   716\\
   145&   213&   324&   416&   567&   625&   713\\
   145&   213&   324&   416&   526&   657&   713\\
   145&   213&   324&   416&   567&   635&   713\\
   145&   213&   324&   416&   536&   657&   713\\
   145&   213&   324&   416&   536&   627&   716\\
   145&   213&   324&   416&   537&   627&   716\\
   145&   213&   324&   416&   526&   637&   716\\
   145&   213&   324&   416&   527&   637&   716\\
   145&   213&   324&   467&   513&   625&   724\\
   145&   213&   324&   467&   567&   635&   724\\
\hline	\end{tabular}\end{table}

\begin{table}[htp]\centering\caption{3+1+1+1+1 Decomposition}\begin{tabular}{|l|l|l|l|l|l|l|}\hline 		
   124&   213&   324&   457&   526&   613&   713\\
   124&   213&   324&   467&   524&   615&   713\\
   124&   213&   324&   467&   524&   615&   715\\
   125&   213&   324&   467&   567&   613&   713\\
   124&   213&   324&   467&   527&   615&   713\\
   124&   213&   324&   467&   524&   635&   715\\
   125&   213&   324&   457&   536&   614&   713\\
   125&   213&   324&   457&   536&   624&   713\\
   125&   213&   324&   456&   537&   637&   714\\
   125&   213&   324&   457&   536&   624&   716\\
   125&   213&   324&   416&   537&   625&   716\\
   125&   213&   324&   416&   547&   625&   716\\
\hline	\end{tabular}\end{table}

\begin{table}[htp]\centering\caption{2+2+1+1+1 Decomposition}\begin{tabular}{|l|l|l|l|l|l|l|}\hline 		
   145&   213&   324&   467&   513&   613&   713\\
   145&   213&   324&   467&   527&   624&   713\\
   156&   213&   324&   467&   527&   624&   713\\
   145&   213&   324&   467&   567&   624&   713\\
   156&   213&   324&   467&   524&   613&   713\\
   145&   213&   324&   467&   513&   625&   713\\
   145&   213&   324&   467&   536&   625&   713\\
   145&   213&   324&   467&   567&   625&   713\\
   145&   213&   324&   467&   526&   635&   713\\
   145&   213&   324&   416&   527&   635&   713\\
   145&   213&   324&   467&   567&   635&   713\\
   145&   213&   324&   467&   513&   625&   725\\
   156&   213&   324&   467&   527&   613&   713\\
\hline	\end{tabular}\end{table}

\begin{table}[htp]\centering\caption{2+1+1+1+1+1 Decomposition}\begin{tabular}{|l|l|l|l|l|l|l|}\hline 		
   145&   213&   324&   467&   526&   613&   713\\
   145&   213&   324&   467&   567&   613&   713\\
   145&   213&   324&   467&   527&   635&   713\\
   156&   213&   324&   456&   527&   627&   713\\
   156&   213&   324&   456&   527&   627&   714\\
   156&   213&   324&   456&   537&   627&   714\\
\hline	\end{tabular}\end{table}

\begin{table}[htp]\centering\caption{1+1+1+1+1+1+1 Decomposition}\begin{tabular}{|l|l|l|l|l|l|l|}\hline 		
 		   145&   213&   345&   427&   526&   613&   713\\
\hline	\end{tabular}\end{table}

\newpage

\begin{center}
\begin{longtable}[htp]{|l|l|l|l|l|l|l|l|}
 \caption{$\Delta$-Reducible $7$-point spaces }\\
\hline $R_1\subset X_7$
&124   &213   &324   &413   &513   &613   &713\\
&124   &213   &324   &413   &513   &613   &724\\
&124   &213   &324   &413   &513   &613   &725\\
&124   &213   &324   &413   &513   &613   &756\\
\hline $R_3\subset X_7$
&124   &213   &324   &413   &513   &624   &716\\
&124   &213   &324   &413   &513   &624   &756\\
\hline $R_5\subset X_7$
&124   &213   &324   &413   &513   &625   &716\\
&124   &213   &324   &413   &513   &625   &725\\
&124   &213   &324   &413   &513   &625   &745\\
&124   &213   &324   &413   &513   &625   &746\\
\hline $R_7\subset X_7$
&124   &213   &324   &413   &516   &625   &713\\
&124   &213   &324   &413   &516   &625   &716\\
&124   &213   &324   &413   &516   &625   &735\\
&124   &213   &324   &413   &516   &625   &736\\
\hline $R_9\subset X_7$
&124   &213   &324   &413   &516   &635   &713\\
&124   &213   &324   &413   &516   &635   &716\\
&124   &213   &324   &413   &516   &635   &724\\
&124   &213   &324   &413   &516   &635   &725\\
\hline $X_4\subset X_7$
&124   &213   &324   &413   &536   &627   &715\\
&124   &213   &324   &413   &567   &615   &713\\
&124   &213   &324   &413   &567   &615   &715\\
&124   &213   &324   &413   &567   &625   &713\\
&124   &213   &324   &413   &567   &625   &715\\
&124   &213   &324   &413   &567   &635   &715\\
\hline $I_3\subset X_7$
&124   &213   &324   &435   &516   &624   &713\\
&124   &213   &324   &435   &516   &624   &716\\
&124   &213   &324   &435   &516   &624   &724\\
&124   &213   &324   &435   &516   &624   &725\\
&124   &213   &324   &435   &516   &624   &735\\
&124   &213   &324   &435   &516   &624   &736\\
\hline $I_{10}\subset X_7$
&124   &213   &324   &435   &516   &625   &713\\
&124   &213   &324   &435   &516   &625   &716\\
&124   &213   &324   &435   &516   &625   &724\\
&124   &213   &324   &435   &516   &625   &725\\
&124   &213   &324   &435   &516   &625   &735\\
&124   &213   &324   &435   &516   &625   &736\\
\hline $I_9\subset X_7$
&124   &213   &324   &435   &516   &635   &713\\
&124   &213   &324   &435   &516   &635   &716\\
&124   &213   &324   &435   &516   &635   &724\\
&124   &213   &324   &435   &516   &635   &725\\
&124   &213   &324   &435   &516   &635   &726\\
&124   &213   &324   &435   &516   &635   &735\\
&124   &213   &324   &435   &516   &635   &746\\
\hline $R_4\subset X_7$
&124   &213   &324   &435   &524   &613   &713\\
&124   &213   &324   &435   &524   &613   &715\\
&124   &213   &324   &435   &524   &613   &726\\
&124   &213   &324   &435   &524   &613   &735\\
&124   &213   &324   &435   &524   &613   &746\\
&124   &213   &324   &435   &524   &613   &756\\
\hline $R_6\subset X_7$
&124   &213   &324   &435   &524   &615   &715\\
&124   &213   &324   &435   &524   &615   &726\\
&124   &213   &324   &435   &524   &615   &736\\
\hline $R_2\subset X_7$
&124   &213   &324   &435   &524   &624   &713\\
&124   &213   &324   &435   &524   &624   &715\\
&124   &213   &324   &435   &524   &624   &716\\
&124   &213   &324   &435   &524   &624   &724\\
&124   &213   &324   &435   &524   &624   &736\\
\hline $X_{5B}\subset X_7$
&124   &213   &324   &435   &524   &627   &716\\
&124   &213   &324   &435   &524   &637   &716\\
&124   &213   &324   &435   &524   &637   &726\\
&124   &213   &324   &435   &524   &647   &716\\
&124   &213   &324   &435   &524   &647   &726\\
&124   &213   &324   &435   &524   &657   &716\\
\hline $I_{12}\subset X_7$
&124   &213   &324   &435   &526   &615   &713\\
&124   &213   &324   &435   &526   &615   &715\\
&124   &213   &324   &435   &526   &615   &724\\
&124   &213   &324   &435   &526   &615   &726\\
&124   &213   &324   &435   &526   &615   &735\\
&124   &213   &324   &435   &526   &615   &736\\
&124   &213   &324   &435   &526   &615   &746\\
\hline $I_5\subset X_7$
&124   &213   &324   &435   &526   &635   &713\\
&124   &213   &324   &435   &526   &635   &715\\
&124   &213   &324   &435   &526   &635   &716\\
&124   &213   &324   &435   &526   &635   &724\\
&124   &213   &324   &435   &526   &635   &726\\
&124   &213   &324   &435   &526   &635   &735\\
&124   &213   &324   &435   &526   &635   &746\\
\hline $I_{13}\subset X_7$
&124   &213   &324   &435   &546   &613   &713\\
&124   &213   &324   &435   &546   &613   &715\\
&124   &213   &324   &435   &546   &613   &724\\
&124   &213   &324   &435   &546   &613   &725\\
&124   &213   &324   &435   &546   &613   &726\\
&124   &213   &324   &435   &546   &613   &735\\
&124   &213   &324   &435   &546   &613   &746\\
\hline $I_8\subset X_7$
&124   &213   &324   &435   &546   &615   &713\\
&124   &213   &324   &435   &546   &615   &715\\
&124   &213   &324   &435   &546   &615   &724\\
&124   &213   &324   &435   &546   &615   &725\\
&124   &213   &324   &435   &546   &615   &726\\
&124   &213   &324   &435   &546   &615   &735\\
&124   &213   &324   &435   &546   &615   &736\\
&124   &213   &324   &435   &546   &615   &746\\
\hline $I_{14}\subset X_7$
&124   &213   &324   &435   &546   &625   &713\\
&124   &213   &324   &435   &546   &625   &715\\
&124   &213   &324   &435   &546   &625   &716\\
&124   &213   &324   &435   &546   &625   &724\\
&124   &213   &324   &435   &546   &625   &725\\
&124   &213   &324   &435   &546   &625   &735\\
&124   &213   &324   &435   &546   &625   &736\\
&124   &213   &324   &435   &546   &625   &746\\
\hline $I_7\subset X_7$
&124   &213   &324   &435   &546   &635   &713\\
&124   &213   &324   &435   &546   &635   &715\\
&124   &213   &324   &435   &546   &635   &716\\
&124   &213   &324   &435   &546   &635   &724\\
&124   &213   &324   &435   &546   &635   &725\\
\hline $I_4\subset X_7$
&124   &213   &324   &456   &513   &613   &713\\
&124   &213   &324   &456   &513   &613   &724\\
&124   &213   &324   &456   &513   &613   &725\\
&124   &213   &324   &456   &513   &613   &756\\
\hline $I_2\subset X_7$
&124   &213   &324   &456   &513   &624   &713\\
&124   &213   &324   &456   &513   &624   &716\\
&124   &213   &324   &456   &513   &624   &724\\
&124   &213   &324   &456   &513   &624   &725\\
&124   &213   &324   &456   &513   &624   &756\\
\hline $I_1\subset X_7$
&124   &213   &324   &456   &524   &624   &713\\
&124   &213   &324   &456   &524   &624   &715\\
&124   &213   &324   &456   &524   &624   &724\\
\hline $R_8\subset X_7$
&125   &213   &324   &435   &514   &613   &713\\
&125   &213   &324   &435   &514   &613   &714\\
&125   &213   &324   &435   &514   &613   &724\\
&125   &213   &324   &435   &514   &613   &726\\
&125   &213   &324   &435   &514   &613   &746\\
\hline $X_{5A}\subset X_7$
&125   &213   &324   &435   &514   &627   &716\\
&125   &213   &324   &435   &514   &637   &716\\
\hline $I_{11}\subset X_7$
&125   &213   &324   &435   &546   &613   &713\\
&125   &213   &324   &435   &546   &613   &714\\
&125   &213   &324   &435   &546   &613   &724\\
&125   &213   &324   &435   &546   &613   &725\\
&125   &213   &324   &435   &546   &613   &726\\
&125   &213   &324   &435   &546   &613   &735\\
&125   &213   &324   &435   &546   &613   &746\\
\hline $I_{15}\subset X_7$
&125   &213   &324   &435   &546   &625   &713\\
&125   &213   &324   &435   &546   &625   &714\\
&125   &213   &324   &435   &546   &625   &716\\
&125   &213   &324   &435   &546   &625   &724\\
&125   &213   &324   &435   &546   &625   &725\\
\hline $I_6\subset X_7$
&125   &213   &324   &436   &536   &625   &713\\
&125   &213   &324   &436   &536   &625   &714\\
&125   &213   &324   &436   &536   &625   &716\\
&125   &213   &324   &436   &536   &625   &725\\
\hline $I_{17}\subset X_7$
&126   &213   &324   &435   &546   &615   &713\\
&126   &213   &324   &435   &546   &615   &714\\
\hline $I_{16}\subset X_7$
&156   &213   &324   &456   &513   &624   &713\\
&156   &213   &324   &456   &513   &624   &714\\
\hline	
\end{longtable}
\end{center}
\end{document}